\documentclass [11pt]{amsart}
\usepackage {amsmath, amssymb, amscd, amsthm, epsfig, color}

\input xy
\xyoption {all}

\def \bb {\mathbb}
\def \s {{\star}}

\def \m {{{\mathcal {M}}(r,d)}}
\def \v {{{\mathcal {V}}}}

\def \e {{{\mathcal {E}}}}
\def \f {{{\mathcal {F}}}}

\def \o {{{\mathcal {O}}}}
\def \mor {{{\text{Mor}}_{d}(C, G(r,N))}}
\def \quot {{{\text{Quot}}_{d} ({\mathcal O}^N, r, C)}}
\def \bp {{{B_d (N, r, C)}}}
\def \morm {{{\text{Mor}}_{d-m}(C, G(r,N))}}
\def \mors {{{\text{Mor}}_{d}(C, G(r,N))^{\text{ss}}}}
\def \moru {{{\text{Mor}}_{d}(C, G(r,N))^{\text{unstable}}}}

\newtheorem{theorem}{Theorem}
\newtheorem {lemma}{Lemma}
\newtheorem {corollary}{Corollary}
\newtheorem {proposition}{Proposition}

\theoremstyle{definition}

\theoremstyle {definition}

\begin{document}
\title [Intersection theory of Quot schemes and moduli of bundles with sections] 
{On the intersection theory of Quot schemes and moduli of bundles with sections} 
\author {Alina Marian}
\address {Department of Mathematics}
\address {Yale University, P.O. Box 208283, New Haven, CT, 06520-8283, USA}
\email {alina.marian@yale.edu}
\date{}

\maketitle

\begin{abstract}

We consider a class of tautological top intersection products on the
moduli space of stable pairs consisting of vector bundles together with
N sections on a smooth complex projective curve C. We show that when N is large,
these intersection numbers can equally be computed on the Grothendieck
Quot scheme of coherent sheaf quotients of the rank N trivial sheaf on
C. The result has applications to the calculation of the intersection
theory of the moduli space of semistable bundles on C.

\end{abstract}

\section{Introduction}

\vspace{0.2in}

We study 
an intersection-theoretic relationship between two moduli spaces associated with bundles on a curve.
We fix throughout a smooth complex projective curve $C$ of genus $g \geq 2$. One of the 
two moduli spaces is the Grothendieck Quot
scheme $\quot$ parametrizing rank $N-r$ degree $d$ coherent sheaf quotients of $\o^N$ on 
$C$. The other is a parameter space, subsequently denoted $B_d (N, r, C)$, for stable
pairs consisting of a rank $r$ degree $d$ vector bundle on $C$ together with $N$ sections. \\

Several moduli spaces of pairs of vector bundles and sections can be constructed gauge-theoretically or by GIT, by varying the notion of 
stability of a pair \cite{bdw}, \cite{thaddeus}. $B_d (N,r, C)$ is one of these, and it is the one most closely related to the moduli space $\m$ 
of semistable rank $r$ and degree $d$ bundles on $C$. Indeed, a point in $B_d (N, r, C)$ consists of a semistable rank $r$ degree $d$ bundle $E$ 
and $N$ 
sections $s_1, \ldots, s_N$ such that any proper subbundle of $E$ whose space of sections contains $s_1, \ldots, s_N$, has slope strictly less 
than the slope of $E$. 
We will assume throughout that $N \geq 2r+1$, and $d$ is as large as convenient. Choosing $d$ large is not restrictive for the moduli space 
$\m$, since $\m$ only depends on $d$ modulo $r$. In case $d> 2r (g-1),$ $\bp$ is 
smooth projective \cite{bdw} and there is a surjective morphism forgetting the sections, $$\pi: B_d (N, r, C) \rightarrow \m.$$

In turn, for $d$ large relative to $N, r, g$, $\quot$ is irreducible, generically smooth, of the expected 
dimension \cite{bdw}. The closed points of $\quot$ are exact sequences 
$ 0 \rightarrow E \rightarrow \o^N \rightarrow F \rightarrow 0$ on $C$ with $F$ coherent of degree $d$ and rank
$N-r$. Thus $\quot$ compactifies the scheme $\mor$ of degree $d$ morphisms from $C$ to the Grassmannian
$G(r, N)$, which is the locus of exact sequences on $C$
with $F$ {\em{locally free}} of degree $d$ and rank $N-r$.\\

The two spaces $\bp$ and $\quot$ are birational, agreeing 
on the open subscheme $U$ of $\quot$ consisting of sequences 
$ 0 \rightarrow E \rightarrow \o^N \rightarrow F \rightarrow 0$
with $E^{\vee}$ semistable: this sequence in $\quot$ corresponds to
the point ${\o^N}^{\vee} \rightarrow E^{\vee}$ in $B_d (N, r, C)$. \\

Both $\quot$ and $\bp$ are fine moduli spaces carrying universal structures. We will denote by 
$$0 \rightarrow \e \rightarrow \o^N \rightarrow \f \rightarrow 0 \, \, \, \text{on} \, \,  \quot \times C$$ 
the universal exact sequence associated with $\quot$ and by
$$\o^N \rightarrow \v \, \, \, \text{on} \, \, \bp \times C.$$ 
the universal bundle with $N$ sections corresponding to $\bp$. Here $\e$ and $\v$ are locally free, and  $\o^N 
\rightarrow \v$ coincides with
${\o^N}^{\vee} \rightarrow \e^{\vee}$ on the common open subscheme $U$ of $\bp$ and $\quot$.  \\

We further fix a symplectic basis $\{1, \delta_j, 1 \leq j \leq 2g, \omega \}$ for the cohomology $H^{\star} (C)$ 
of the curve $C$ and let 
\begin{equation}
\label{kun1}
c_i(\e^{\vee}) = a_i \otimes 1 + \sum_{j =1}^{2g} b_i^j \otimes \delta_j + f_i \otimes \omega, \, \, 1\leq i \leq r, \, 
\, \, \text{and}
\end{equation}  
\begin{equation}
\label{kun2}
c_i(\v) = {\bar{a}}_i \otimes 1 + \sum_{j =1}^{2g} {\bar{b}}_i^j \otimes \delta_j + {\bar{f}}_i \otimes 
\omega, \, \, 1\leq i \leq r,
\end{equation}
be the K\"{u}nneth decompositions of the Chern classes of $\e^{\vee}$ and $\v$ with 
respect to this basis
for $H^{\star} (C).$ Note that $$a_i \in H^{2i} (\quot), \, b_i^j \in H^{2i-1}(\quot), \, f_i \in H^{2i-2} (\quot),$$
and ${\bar{a}}_i, {\bar{b}}_i^j, {\bar{f}}_i$ belong to the analogous cohomology groups of $\bp.$ \\

The goal of this paper is to present the following observation.

\begin{theorem}
Let $r\geq 2$ and $N\geq 2r+1.$ Assume moreover that $C$ has genus at least 2 and that $d$ is large enough so that both $\quot$ and $\bp$ are 
irreducible of the expected 
dimension 
$ e= Nd -r (N-r) (g-1).$ Let $s$ and $t$ be 
nonnegative integers such that $t \leq \frac{N}{r}$ and s r + t = $e$. Let $P(a, f, b)$ be a polynomial of degree 
$t$ in the 
$a$, $f$ and (algebraic combinations) of the $b$ classes. Then 
\begin{equation}
\label{transfer}
\int_{\bp} {\bar{a}}_r^s P(\bar{a}, \bar{f}, \bar{b} ) = \int_{\quot} a_r^s P(a, f, b).
\end{equation}
\end{theorem}

\vspace{0.1in}

Note that the theorem by no means accomplishes a complete comparison of the intersection theories
of $\bp$ and of $\quot$. Its scope is in fact modest: it asserts that the intersection of any
algebraic polynomial $P(a, f, b)$ in the tautological classes paired with a complementary self-intersection
$a_r^s$ can be evaluated to the same result on $\quot$ and on $\bp$ only asymptotically as $N$ is large 
relative to the degree of $P$. A more general statement is false due to the different nature of the nonoverlapping loci of
$\quot$ and $\bp$. \\

The proof of the theorem, which will occupy the next section, relies on an explicit realization of the tautological
classes on $\bp$ and $\quot$ as degeneracy loci (or pushforwards of degeneracy loci) of morphisms stemming from the 
universal $\o^N \rightarrow \v$ on $\bp \times C$ and $\o^N \rightarrow \e^{\vee}$ on $\quot \times C.$ 
This approach was used in \cite{bertram} to study the intersection theory of the tautological 
$a$ classes on $\quot$. Bertram proved that generic degeneracy-loci representatives of the $a$ classes in a top-degree
polynomial $P(a_1, \ldots, a_r)$ intersect properly along $\mor$ and avoid the boundary $\quot \setminus \mor$. Thus the 
intersection $P(a_1, \ldots, a_r) \cap \quot$ has enumerative meaning, and expresses a count of maps from $C$ to $G(r,N)$ 
satisfying incidence properties. \\

Similarly, in the situation described by Theorem 1, we will show that generic degeneracy-loci representatives of the 
$a$, $f$, $b$ classes in the product 
$a_r^s P(a, f, b)$ on the one hand intersect properly along the common open subscheme $U$ of $\quot$ and $\bp$ and 
on the other hand manage to avoid the 
nonoverlapping loci $\bp \setminus U$ and \\ 
$\quot \setminus U.$ The argument uses heavily the methods and results of \cite{bertram}. \\

The theorem was sought, and is indeed effective, in the attempt to calculate the intersection 
theory of the moduli space $\m$ of semistable bundles on $C$, by transferring intersections to
the Quot scheme $\quot$ and computing them there. This idea was inspired by \cite{witten}.
The intersection-theoretic relationship between $\m$ and $\quot$ is simplest when the rank $r$ and 
the degree $d$ are coprime. In this case, 
$\m$ is a smooth projective variety and a fine moduli space with a universal bundle $V \rightarrow \m \times 
C.$  Furthermore, $\bp$ is a projective bundle over $\m$,
$$\pi: \bp = {\bb P} ({\eta_{\star} V}^{\oplus N}) \rightarrow \m ,$$  where $\eta: \m \times C \rightarrow 
\m $ is the projection. Intersections on $\m$ can therefore be easily lifted to $\bp.$ To be
precise, let  
$$c_i(V) = {\tilde{a}}_i \otimes 1 + \sum_{j = 1}^{2g} {\tilde{b}}_i^j \otimes \delta_j + {\tilde{f}}_i \otimes \omega, \, \, 1\leq i
\leq r,$$ be the K\"{u}nneth decomposition of the Chern classes of $V,$ and let $$\zeta = c_1 (\o (1)) \, \, \,
\text{on} \, \, \bp = {\bb P} ({\eta_{\star} V}^{\oplus N}).$$ It is easy to see that $\v = \pi^{\star} V \otimes \o(1)$ on $\bp \times C.$ We 
choose $N$ so that the 
dimension $N (d-r(g-1)) -1$ of the fibers of $\pi$ is a multiple of $r$, $N(d -r (g-1)) -1 = rs.$ Then if $P(\tilde{a}, \tilde{b}, 
\tilde{f})$ is a top-degree 
polynomial on $\m$ invariant under 
twisting the universal $V$ with line bundles from $\m$, we have 
\begin{equation}
\int_{\m} P(\tilde{a}, \tilde{b}, \tilde{f}) = \int_{\bp} \zeta^{\text{dim fiber}}P(\bar{a}, \bar{b}, \bar{f}) = \int_{\bp} 
{\bar{a}}_r^s
 P(\bar{a}, \bar{b}, \bar{f}).
\end{equation} 
We therefore obtain the following consequence of Theorem 1.
\begin{corollary}
For $r$ and $d$ coprime, top intersections on $\m$ can be viewed as intersections on $\quot$ in the regime
when $N$ is large relative to the dimension of $\m$, and $d$ is large enough so that $\quot$ is irreducible. 
Specifically, for $N(d -r (g-1)) -1$ 
divisible by $r$, $N \geq r \cdot \dim \m,$ and $d \geq (2g-1) (N+2)^{r-1},$ we have 
\begin{equation}
\int_{\m} P(\tilde{a}, \tilde{b}, \tilde{f}) = \int_{\quot} a_r^s  P(a, b, f).
\end{equation}
\end{corollary} 

The calculation of the tautological intersection theory on 
the Quot scheme is manageable
via equivariant localization with respect to a natural ${\bb C}^{\star}$ action \cite{quot}. Theorem 1 together with the results of 
\cite{quot} thus constitute an effective method for evaluating the top pairings of the cohomology generators of $\m$, and provide an 
alternative 
to the important work \cite{jeffkir}, \cite{kefeng} on the intersection theory of $\m$. 
In the rank 2, odd degree  
case, complete calculations were carried out in \cite{om}. A derivation of the volume of $\m$ as well as of other intersection numbers for 
arbitrary $r$ is pursued in \cite{volumes}. \\

\vspace{0.3mm}

{\bf {Acknowledgements.}} The research for this paper was carried
out during the author's 
graduate work at Harvard and she would like to thank her advisor, 
Shing-Tung Yau, for his guidance and support. Conversations with Dragos Oprea were as usual pleasant and instructive.\\

\section{Tautological intersections on $\quot$ and $\bp$}

\vspace{0.2in}

Following \cite{bertram} closely, we first summarize the main facts concerning the
boundary of the compactification by $\quot$ of the space $\mor$ of maps from $C$ to the Grassmannian $G(r,N)$. We
then analyze the relevant intersections on $\quot$ and $\bp$, with conclusions gathered in Proposition 1 and Proposition 2
respectively. Theorem 1 follows from the two propositions. \\

\subsection{The Quot scheme compactification of the space of maps to the Grassmannian $G(r,N)$} $\quot$ parametrizes exact sequences $0 \rightarrow E
\rightarrow \o^N \rightarrow F \rightarrow 0$ with $F$ coherent of degree $d$ and rank $N-r$. Two quotients $\o^N \rightarrow F_1 \rightarrow 0$ and
$\o^N \rightarrow F_2 \rightarrow 0$ are equivalent if their kernels inside $\o^N$ coincide. Now let $0 \rightarrow S \rightarrow \o^N \rightarrow Q
\rightarrow 0$ denote the tautological sequence on $G(r,N)$. By pullback of the tautological sequence to $C$, a degree $d$ map $f: C \rightarrow
G(r,N)$ gives a short exact sequence $0 \rightarrow E \rightarrow \o^N \rightarrow F \rightarrow 0$ of vector bundles on $C$ and
conversely.

Thus $\mor$ is contained in $\quot$ as the locus of {\em{locally free}} quotients 
of $\o^N$ on $C$. \\

There are natural evaluation morphisms
\begin{equation}
ev: \mor \times C \rightarrow G(r,N), \, \, \, (f,p) \mapsto f(p),
\label{ev}
\end{equation}
 and for $p \in
C,$
\begin{equation}
ev_{p}: \mor \rightarrow G(r,N), \, \, \, f \mapsto f(p).
\label{evp}
\end{equation}

The Zariski
tangent space at a closed point $f \in \mor$ is $H^0 (C, f^{\s} TG(r, N)).$
The expected
dimension, which we subsequently denote by $e$, is therefore
$$e = \chi(C, f^{\s} TG(r, N))=
N d -r(N-r) (g-1),$$
by the Riemann-Roch formula on $C$. In {\cite{bdw}},
Bertram,
Daskalopoulos and Wentworth established\\

{\it For sufficiently large degree $d$
relative to $r$, $N$, and the genus $g$ of $C$, $\mor$ is an irreducible quasiprojective
scheme of dimension as expected, $e$. In this regime, the compactification $\quot$ is an irreducible
generically smooth projective scheme.} \\

We assume throughout this paper that the degree is very large compared to $r$, $N$ and $g$. \\

What remains in
$\quot$ once we take out $\mor$ is the locus of quotients of $\o^N$ which are not
locally free. Since a torsion-free coherent sheaf on a smooth curve is locally free,
these nonvector bundle quotients are locally free modulo a torsion subsheaf which they
necessarily contain. That is to say, for each such $\o^N \rightarrow F \rightarrow 0$
there
is an exact sequence $$0 \rightarrow T \rightarrow
F \rightarrow F/T \rightarrow 0$$
with $F/T$ locally free of degree at least 0, and $T$ supported at finitely many points. 
The degree of $T$ can therefore be at most
$d$. Let
$B_m$ be the locus in $\quot$ of quotients whose torsion subsheaf $T$ has degree $m$.
$$\quot = \mor \sqcup_{m = 1}^{d} B_m,$$
with the $B_m$s playing the role of {\em{boundary}} strata. 
We will denote by ${\mathcal B}$ the
union $\sqcup_{m = 1}^{d} B_m.$ \\

We analyze ${\mathcal B}$ more carefully.
If $F$ is a quotient with torsion $T$ of degree $m$, $F/T$ is a 
vector bundle of rank $N-r$ and degree $d-m.$ By associating to 
$\o^N \rightarrow F \rightarrow 0$ the quotient $\o^N \rightarrow 
F/T \rightarrow 0$, we get a morphism 
\begin{equation}
\label{P}
P: B_m \rightarrow \morm, \, \, \, \{\o^N \rightarrow F \}\, \longmapsto \, \{\o^N \rightarrow
F/T \}.
\end{equation}
Let us consider the
inverse image under $P$ of an element $$0 \rightarrow E' \rightarrow \o^N \rightarrow
F' \rightarrow 0$$ in $\morm$. 
For each
sheaf morphism $E' \rightarrow T$, with $T$ torsion of degree $m$, there is a sheaf
$F$ constructed as the pushout of the maps $E' \rightarrow \o^N$ and $E' \rightarrow T$
in the category of coherent sheaves on $C$. Thus the
following diagram of sheaves can be completed, if we start with the data of
the first row and of the first column map:
$$
\begin{array}{ccccccccc}
0 & \longrightarrow & E' & \longrightarrow & {\mathcal {O}}^N & \longrightarrow &   
F' & \longrightarrow & 0 \\
 & & \downarrow & & \downarrow & & \| & & \\
0 & \longrightarrow & T & \longrightarrow & F & \longrightarrow & F' &
\longrightarrow & 0
\end{array}
$$
Moreover the map $\o^N \rightarrow F$ is surjective if $E' \rightarrow T$ is. If we
denote its kernel by $E$ we have

$$
\begin{array}{ccccccccc}
 & & E  & = & E & & & & \\
 & & \downarrow & & \downarrow & & & & \\
0 & \longrightarrow & E' & \longrightarrow & {\mathcal {O}}^N & \longrightarrow &
F' & \longrightarrow & 0 \\
 & & \downarrow & & \downarrow & & \| & & \\
0 & \longrightarrow & T & \longrightarrow & F & \longrightarrow & F' &
\longrightarrow & 0
\end{array},
$$
with the first two columns also being short exact sequences.

Finally assume that two quotients $\o^N \rightarrow F_1 $ and $\o^N \rightarrow F_2$ induce
quotients $\o^N \rightarrow F_1/T_1$ and $\o^N \rightarrow F_2/T_2$ which have the same 
kernel $E'$ inside $\o^N.$ It is then immediate that 
%It is also immediate that two degree $d$ quotients $\o^N \rightarrow
%$ F_1$ and
%$\o^N \rightarrow F_2$, such that the induced quotients $$\o^N \rightarrow F_1/T_1$$ and
%$$\o^N \rightarrow F_2/T_2$$ have the same kernel $E'$, represent the same element of
%$\quot$ if and only if the map $E' \rightarrow T_2$ is obtained from $E' \rightarrow
%T_1$ by composing it with an isomorphism $T_1 \rightarrow T_2$. That is to say,
$\o^N
\rightarrow F_1$ and $\o^N \rightarrow F_2$ represent the same element of $\quot$ if
and only if $E' \rightarrow T_1$ and $E' \rightarrow T_2$ represent the same element of
the Quot scheme of degree $m$, rank 0 quotients of $E'$. This latter scheme, which we
denote by ${\text{Quot}}_{m} (E', r, C)$, is a
smooth projective variety of dimension $rm$. \\

We conclude that the morphism
$$P: B_m \rightarrow \morm$$ is surjective with fiber over
 $$0 \rightarrow E' \rightarrow
\o^N \rightarrow F' \rightarrow 0$$ isomorphic to ${\text{Quot}}_{m} (E', r, C).$
From this it follows easily that the codimension of the boundary stratum $B_m$ in
$\quot$, for $m$
small, is equal to 
\begin{equation}
\label{codim0}
\text{codim}\, B_m = (N-r)m.
\end{equation}
The codimension of $B_m$ for $m$ large (so that $\morm$ is not
irreducible) increases linearly with the degree $d$.\\

Note also that by associating to the
quotient $\o^N \rightarrow
F$ the {\em{double}} data of
$\o^N \rightarrow F/T$ {\em{and}} of the support of $T$, we get a morphism
\begin{equation}
\label{Pl}
Pl: B_m \rightarrow \text{Sym}^m\, C  \times \morm,
\end{equation}
where $\text{Sym}^m \, C$ denotes the $m$th symmetric product
of the curve $C$. \\

\subsection{Intersections on the Quot scheme.}

The classes that we would like to intersect are the K\"{u}nneth components \eqref{kun1} of $c_i (\e^{\vee}), 1\leq i \leq r.$ Let
$$\eta: \quot \times C \rightarrow \quot$$ be the projection, and note that
$$a_i = c_i (\e_p^{\vee}), \, \, \text{for} \, \, p \in C, \, \, f_i = \eta_{\star} (c_i(\e^{\vee})), \, \, 1\leq i \leq r.$$ 
Furthermore, combinations of 
$b$ classes which are algebraic can be expressed in terms of $a$ classes and of 
pushforwards under $\eta$ of polynomials in the Chern classes of $\e^{\vee}$.
For instance,
$$2\sum_{j=1}^{g} b_1^j b_1^{j+g} = 2d a_1  - \eta_{\star}(c_1^2),$$

$$ \sum_{j=1}^{g} b_1^j b_2^{j+g} + \sum_{j=1}^{g} b_2^j b_1^{j+g} = a_1 \eta_{\star} (c_2) + d a_2  - \eta_{\star} (c_1 c_2).$$

We will view an arbitrary homomorphism
$$P(a,b, f): A_{\star} (\quot) \rightarrow A_{\star - \deg P} (\quot)$$ as an actual intersection operation 
of subschemes. For simplicity we will omit $b$ classes from the discussion, since modulo $a$ classes, they are pushforwards
of polynomials in $c_i(\e^{\vee})$, hence can be treated quite similarly to the $f$ classes.  \\

{\bf Notation.} For each $i, 1\leq i \leq r$, we choose $\o^{r-i+1}$, a trivial subbundle of $\o^N$ on $\quot \times C.$ 
The universal 
exact sequence $0 \rightarrow \e \rightarrow \o^N \rightarrow \f \rightarrow 0$ on $\quot \times C$ induces for each 
$i$, by dualizing and restricting, morphisms 
$$\varphi_i: \, \o^{r-i+1} \rightarrow \e^{\vee}.$$ We let $D_i$ be the subscheme of $\quot \times C$ which is the 
noninjectivity locus of $\varphi_i$. Similarly, by fixing $p \in C$, we denote by $D_{i,p}$ the subscheme of $\quot$ 
which is the degeneracy locus of the restriction $\varphi_{i,p}$ of $\varphi_i$ to $\quot \times p.$ We let ${\bb 
D}_i \in
A_{e-i} (D_i)$ and ${\bb D}_{i,p} \in A_{e-i} (D_{i,p})$ be the $i$th degeneracy {\em {classes}} of $\varphi_i$ and
$\varphi_{i,p}$ respectively, satisfying
$$ {\bb D}_i = c_i (\e^{\vee}) \cap [\quot \times C], $$ $${\bb D}_{i,p} = c_i (\e_p^{\vee}) \cap [\quot].$$
Denote by $D_i^0$ and $D_{i,p}^{0}$ the restrictions of the degeneracy loci to $\mor \times C$ and $\mor$ 
respectively. Further, we let $Y_i$ be the degeneracy locus of the map
$$\o^{r-i+1} \rightarrow S^{\vee} \, \, \, \text{on} \, \, \, G(r,N).$$

It is well known that $Y_i$ represents the $i$th Chern class of $S^{\vee}$, and therefore has codimension $i$ in $G(r,N).$ 
Now $D_i^0$ and $D_{i,p}^0$ are precisely the pullbacks of $Y_i$ by the evaluation maps \eqref{ev}, \eqref{evp},
\begin{equation}
D_i^0 = ev^{-1}(Y_i), \, \, \, \, \, \, D_{i,p}^0 = ev_{p}^{-1} (Y_i). 
\end{equation}
Thus $D_i^0$ and $D_{i,p}^0$ have codimension $i$ in $\mor \times C$ and $\mor$ respectively, and $D_{i,p}^0$ is the locus of degree $d$ maps from 
$C$ to $G(r,N)$ which send the point $p$ to $Y_i$. Likewise, $\eta(D_i^0)$ is the locus of maps whose images intersect $Y_i,$ 
and $\eta (D_i^0)$ has codimension $i-1$ in $\mor$. 
This last statement follows because the fibers of
the restricted $\eta: D^0_{i} \rightarrow \eta (D^0_{i})$ consist of finitely
many
points, outside a locus of high codimension in $\eta (D^0_{i})$. (This is
the
locus of maps $f: C \rightarrow G(r,N)$ such that $f(C)$ is contained in $Y_{i},$
that is to say the morphism scheme ${\text{Mor}}_{d} (C, Y_{i})$; its codimension
in $\mor$ is proportional to the degree $d$ of the maps, which is assumed very large.) 
In fact, for $i \geq 2,$ the fibers of the morphism $\eta: D_i^0 \rightarrow \eta 
(D_i^0)$
consist generically of exactly one point of the curve $C$. To see this,
consider the `double evaluation' map $$ev^2: \mor \times C \times C \rightarrow
G(r,N) \times G(r,N), \, \, \, (f, p, q) \longmapsto (f(p), f(q)).$$ Let $\eta^2$ be
the projection $$\mor \times C \times C \rightarrow \mor.$$ The subvariety of
the base $\eta (D^0_{i})$ over which $\eta$ has fibers consisting of at least two
points is $\eta^2 ({\left (ev^2 \right )}^{-1} (Y_{i} \times Y_{i}))$ and has codimension $2i-2$ in
$\mor$, thus codimension $i-1$ in $\eta (D^0_{i})$ (here $i \geq 2.$) The fibers of
$\eta: D^0_{i} \rightarrow \eta (D^0_{i})$ consist therefore generically of 
exactly one point, and $\eta(D_i^0)$ has codimension $i-1$ in $\mor$. \\

The codimension \eqref{codim0} of the boundary stratum $B_m$ is $(N-r) m$ and $(N-r)m  > i$, for $1 \leq i \leq r$, $1\leq
m \leq d$ and
for $N \geq 2r
+1.$ As the codimension of each irreducible component of the degeneracy locus $D_i$ is at most $i$ on general grounds \cite{fulton}, Theorem 14.4, 
we conclude\\

{\it $D_i$  has codimension $i$ in $\quot \times C$, $D_{i,p}$ has codimension $i$ in $\quot$, and $\eta(D_i)$ has
codimension
$i-1$ in $\quot$.}\\

Now $GL(N)$ acts on $\quot$ and on $\mor$ as the automorphism group of $\o^N$:
for $g \in
GL(N),$ $$[\pi: \o^N \rightarrow F ] \mapsto [\pi \circ g:
\o^N \rightarrow F].$$ 
Translating degeneracy loci by the action of the group enables us to make intersections on $\mor$ proper. Indeed, for any 
subvariety $Z$ of $\mor$ and a general
$g \in GL(N)$, the intersection $Z \cap g D^0_{i, p}$ is either empty or has  
codimension $i$ in $Z$: by Kleiman's theorem \cite{kleiman}, for a general $g$, the
intersection $ev_{p} (Z) \cap g Y_{i}$ is either empty or of the expected dimension
in the Grassmannnian $G(r,N)$, whence the intersection $Z \cap ev_{p}^{-1}(g
Y_{i})$ also has the expected dimension. Similarly, for
a general $g \in GL(N)$ the
intersection $Z \times C
\cap ev^{-1} (g Y_{i})$ is either empty or
 has codimension $i$ in $Z \times C.$ Since for a general $g$ the generic fiber of
the projection $$\eta: Z \times C \cap
ev^{-1} (g
Y_{r-i+1}) \rightarrow Z \cap \eta ev^{-1} (g Y_{i})$$ is finite,
the intersection
$Z \cap \eta ev^{-1} (g Y_{i}) = Z \cap \eta (g D_i^0) $ has codimension $i-1$ in $Z$.

We let $\mors$ denote the smooth open subscheme of $\mor$ whose closed points are
exact sequences $0 \rightarrow E \rightarrow \o^N \rightarrow F \rightarrow 0$ with
$E^{\vee}$ semistable. The argument above shows
that for a generic choice of trivial subbundle $\o^{r-i+1}$, $D_{i,p}^0$ intersects the unstable locus 
$$\moru = \mor \setminus 
\mors$$ properly. For dimensional reasons, the restriction map
$$A_{e-i}(D_{i,p}) \longrightarrow A_{e-i} (D_{i,p}^0 \cap \mors )$$ is therefore an isomorphism.
We conclude that the degeneracy {\em {class}} ${\bb D}_{i,p}$ is not only in the top Chow group 
of $D_{i,p}$, but in 
fact
$${\bb D}_{i,p} = [D_{i,p}]$$ since the equality holds for the restrictions to the smooth open subscheme $\mors$ of 
$\mor$.  

Similarly,  $${\bb D}_{i} = [D_{i}], $$ and as $\eta: D_i \rightarrow \eta (D_i)$ is generically one-to-one, we have
$$\eta_{\star}[D_i] = [\eta (D_i)].$$ 

We collect these conclusions in the statement of the following lemma.
\begin{lemma}
(1) For a generic choice of trivial subbundle $\o^{r-i+1}$ of $\o^N$ on $\quot \times C$, 
$$[D_i] = c_i (\e^{\vee}) \cap [\quot \times C ],$$ $$[D_{i,p}] = a_i \cap [\quot]\, \,
\, \text{and}$$
$$[\eta (D_i)] = f_i \cap [\quot].$$
(2) Generic representative subschemes for all instances of $a_i$, $f_i$ and algebraic combinations of $b_i^j$ classes in a monomial $P(a, b, f)$ 
intersect the unstable locus of $\mor$ and each other properly along $\mor$. 
\label{deg22}
\end{lemma}

We would now like to show that the intersection of generic degeneracy loci representatives of the
tautological classes appearing in a monomial $P (a, b, f) a_r^s,$ for $e = sr + t$, $t= \deg P$ and $t \leq \frac{N}{r}$
avoids the boundary $\sqcup_{m = 1}^{d} B_m$ of the Quot scheme. Lemma 1 ensures on the other hand that intersections occur with
the expected codimension along $\mor$, even avoiding the unstable locus of $\mor$. We will thus be able to conclude that the intersection
number $P(a,b,f) a_r^s \cap [\quot ]$ is the number of intersection points in $\mors$ (counted with multiplicities) of generic 
degeneracy loci representatives of all the classes in $P(a,b,f) a_r^s$. \\

We therefore analyze the intersection of a degeneracy locus 
$D_{i,p}$ with a boundary stratum $B_m$. The 
subtlety here is that although the $B_m$s have large codimension in $\quot$, subschemes of the type $D_{i,p}$ and $\eta (D_i)$ may not intersect them 
properly, 
therefore giving an overall excess intersection with the boundary. 

Suppose $0 \rightarrow E \rightarrow \o^N \rightarrow F \rightarrow 0$ is a point in the intersection 
$D_{i,p} \cap B_m.$ This means that $F$ has a torsion subsheaf $T$ of degree $m$ and the vector bundle map 
$\o^{r-i+1} \rightarrow E^{\vee}$ is not injective at $p$. Recall the surjective morphism
\eqref{Pl},
$$Pl: B_m \rightarrow \text{Sym}^m\, C  \times \morm, \, \, \, \{\o^N \rightarrow F\} \longmapsto \{ \text{Supp}\, T , \, \o^N 
\rightarrow F/T \},$$
and its further composition $P$ with the projection to $\morm$, given by \eqref{P}.  
Let $$0 \rightarrow \e^m \rightarrow \o^N \rightarrow \f^m \rightarrow 0$$ be the universal 
sequence on $\morm \times C$ and denote by $D_{i,p} (m)$ the degeneracy locus of $\o_p^{r-i+1} \rightarrow {\e_p^m}^{\vee}$ on 
$\morm$. 

Let $E'$ be the kernel of the map $\o^N \rightarrow F/T.$ If $p \notin \text{Supp} \, T$, then the map $\o^{r-i+1} \rightarrow {E'}^{\vee}$
is also not injective at $p$ since it coincides with $\o^{r-i+1} \rightarrow E^{\vee}$ outside the support of $T$. We conclude that
\begin{equation}
\label{boundaryint}
D_{i,p} \cap B_m \subset P^{-1} (D_{i,p} (m)) \cup Pl^{-1} (p \times \morm).
\end{equation}

We pick now $s$ distinct points $p_k$ and $s$ generic $GL(N)$ elements $g_k$, $1\leq k \leq s.$ Each instance of
$a_r$ in the product $a_r^s$ is represented by a subscheme 
$g_k D_{r, p_k}, 1\leq k \leq s.$ 
We want to show that 
\begin{equation}
\label{emptyset}
\cap_{k = 1}^{s} g_k D_{r, p_k}
\cap B_m = \emptyset.
\end{equation}
Indeed, if we let $I$ index subsets of $\{1, \ldots, s\}$ of cardinality $s-m$, then \eqref{boundaryint} implies
that 
\begin{equation}
\cap_{k = 1}^{s} g_k D_{r, p_k} \cap B_m \subset \bigcup_{I} P^{-1}(\cap_{k \in I} g_k
D_{r, p_k} (m)).
\end{equation}
For general $g_k$s, the subschemes $$g_k
D_{r, p_k} (m), \, 1 \leq k \leq s$$ intersect properly on
$\text{Mor}_{d-m}.$ The codimension of $\cap_{k \in I} g_k D_{r, p_k} (m)$ in $\morm$ is therefore 
\begin{equation}
\label{codimension}
(s-m) r = e - t - mr \geq N d - r(N-r) (g-1) -\frac{N}{r} - mr.
\end{equation}
The inequality above holds in the regime that we are interested in, 
when $e - rs = t \leq \frac{N}{r}.$ 

A dimension count for $\morm$ will show now that 
$\cap_{k \in I} g_k D_{r, p_k} (m)$ is empty for all $I$ with $|I| = s-m$. 

There are two cases to be
considered.\\

{\em{Case 1.}} $m$ is large enough so that $\text{Mor}_{d-m}$ is not irreducible of
the expected dimension. The dimension of $\morm$ can be bounded in this case independently of 
the degree $d$. The codimension \eqref{codimension} of  $\cap_{k \in I} g_k D_{r, p_k} (m)$
is however bounded below by an expression increasing linearly with $d$.
For $d$ sufficiently large 
therefore, the intersection $\cap_{k\in I } g_k D_{r, p_k} (m)$ is empty.\\

{\em{Case 2.}} $m$ is so that $\text{Mor}_{d-m}$ is irreducible of the expected
dimension $$N(d-m) -r (N-r) (g-1).$$
Comparing this dimension with the codimension \eqref{codimension} of
$\cap_{k \in I } g_k D_{r, p_k} (m)$, we find
$$N(d-m) -r (N-r) (g-1) < N d - r(N-r) (g-1) -\frac{N}{r} - mr, \, \forall 1\leq m \leq d, \, \, \text{if} \, N \geq 2r+1. $$
In this case also therefore, the intersection $\cap_{k=1}^{s} g_k D_{r, p_k} \cap B_m$ is empty.\\

Taking into account the analysis leading to Lemma 1, we have the following statement. 

\begin{proposition} The intersection of generic degeneracy loci representatives for the classes appearing in the top product $P(a,f,b) a_r^s, \, \deg P \leq 
\frac{N}{r},$ on $\quot$ avoids the boundary ${\mathcal B} = \sqcup_{m = 1}^{d} B_m$ of $\quot$. Furthermore, such degeneracy loci meet each 
other as well as $\moru$ properly along $\mor$. Their scheme-theoretic intersection is therefore zero-dimensional of length equal to the evaluation
$P(a,f,b) a_r^s \cap [\quot]$, and lies entirely in $\mors.$
\end{proposition}

\subsection{Intersections on the moduli space of stable pairs.} For $d > 2r (g-1),$ the moduli space $\bp$ is a 
smooth projective variety which parametrizes morphisms $\varphi: \o^N \rightarrow E$ on $C$ with $E$ being a semistable
rank $r$ and degree $d$ bundle, and with the property that the subsheaf ${\text{Im}}\, \varphi$ of $E$ is not destabilizing. In other words,
if ${\text{Im}}\, \varphi$ has degree $d'$ and rank $r'$, there is a strict inequality
\begin{equation}
\label{strictin}
\frac{d'}{r'} < \frac{d}{r}.
\end{equation}
This moduli space was described, together with other related moduli spaces of pairs of bundles and sections, in \cite{bradlow},\cite{bdw}.
The case of bundles together with {\em{one}} section was studied in \cite{thaddeus}.  \\

$\bp$ is a fine moduli space carrying a universal morphism
$${\Phi} : \o^N \rightarrow \v \, \, \, \text{on}\, \, \, \bp \times C.$$

We let $U$ be the open subscheme of $\bp$ consisting of generically surjective (as maps of
vector bundles) morphisms $\o^N \rightarrow E.$ 
Crucially,
$U$ coincides with
the semistable locus of $\quot$, {\em{i.e.}} with the open subscheme consisting
of exact
sequences $0 \rightarrow E \rightarrow \o^N \rightarrow F \rightarrow 0$ with semistable
$E^{\vee}.$ The universal structures on
$\quot$ and $\bp$ agree on $U$: $${\o^N}^{\vee} \rightarrow \e^{\vee}
\text{ on }\; U \times
C \; \text{coincides with} \; \o^N \rightarrow \v \; \text{on}\; {U}
 \times C.$$  The actions of $GL(N)$
on $\quot$ and $\bp$ also naturally coincide on $U$. \\

In order to study some of the intersection theory of $\bp$ in relation to that
of the Quot scheme, it will be necessary to understand the boundary $\bp \setminus U.$ Let 
$B(r', d'), \, \, r' < r$ be the locus of maps $\varphi: \o^N \rightarrow E$ such that the locally free 
subsheaf ${\text{Im}}\, \varphi$ of $E$ has rank $r'< r$ and degree $d'$. There are finitely many such loci since
$$0 \leq d' <  \frac{r'}{r} d \, \, \text{and} \, \, \,  1\leq r' \leq r-1.$$ Moreover, 
\begin{equation}
\label{flat}
\bp = U \sqcup_{r'< r, d'} B(r', d').
\end{equation} 
The $B(r',d')$s for $r' < r$ are strata in the flattening stratification
of $\bp$ for the sheaf ${\text{Im}}\, \Phi$ and the projection (slightly abusing notation) 
$\eta: \bp \times C \rightarrow \bp.$  \\

We aim to investigate the codimension of $B(r', d')$ in $\bp$. A quotient $\varphi: \o^N \rightarrow E$ in $B (r', d')$ gives on the one hand a 
point 
${{\text{Im}}\, \varphi}^{\vee} \hookrightarrow {\o^N}^{\vee}$ in the Quot scheme
${\text{Quot}}_{d'} (\o^N, r', C)$ of rank $r'$, degree $-d'$ subsheaves of $\o^N$. We denote this morphism
by $\Psi_{r', d'},$ 
\begin{equation}
\label{psi}
\Psi_{r', d'}: B(r', d') \rightarrow {\text{Quot}}_{d'} (\o^N, r', C), \, \, \, \{\varphi: \o^N \rightarrow E \} \mapsto 
\{{{\text{Im}}\, \varphi}^{\vee} \hookrightarrow {\o^N}^{\vee}\}.
\end{equation} 
On the other hand, $\varphi: \o^N 
\rightarrow E$ determines a 
point ${\text{Im}}\, \varphi \hookrightarrow E$ in the Quot 
scheme
${\text{Quot}}_{d-d'}(E, r', C)$ of rank $r'$, degree $d'$ subsheaves of $E$. Conversely, a point in $B(r', d')$ is determined by a choice of 
semistable $E$ in the moduli space
$\m$ of semistable rank $r$ and 
degree $d$ bundles on $C$, and a choice of points in ${\text{Quot}}_{d-d'}(E, r', C)$ and ${\text{Quot}}_{d'} (\o^N, r', C)$ such that the subsheaf of
$E$ is isomorphic to the dual of the subsheaf of $\o^N.$
Now \cite{mihnea} provides an upper bound for the dimension of ${\text{Quot}}_{d-d'}(E, r', C)$ which is independent of the choice
of $E$ inside $\m$, 
$$\dim {\text{Quot}}_{d-d'} (E, r', C) \leq r'(r-r') + d r' - d' r.$$ 
An upper bound for the dimension of $B(r', d')$ is then 
$$\dim B(r', d') \leq \dim {\text{Quot}}_{d'} (\o^N, r', C) + r'(r-r') + dr'-d'r + \dim \m.$$ 
Recall that $\dim \bp = N d - r(N-r)(g-1),$ and as established in \cite{bdw}, $\dim {\text{Quot}}_{d'} (\o^N, r', C)$  either has a 
degree-independent upper 
bound (for low $d'$)
or else $$\dim {\text{Quot}}_{d'} (\o^N, r', C) = N d' -r' (N-r') (g-1).$$ 
A lower bound for the codimension of $B(r', d')$ in $\bp$ is thus 
$$
{\text{codim}} \,B(r', d') \geq N d - N d' - d r' +d' r + \, \, \text{degree-independent terms}.
$$
When $d$ is very large, it follows that 
\begin{equation}
\label{codimd}
{\text{codim}} \,B(r', d') \geq  (N-r) (d-d') > \frac{N-r}{r} d.
\end{equation}
For the last inequality we took account of the stability condition for morphisms $\varphi: \o^N \rightarrow E$ 
expressed by \eqref{strictin}.\\

We conclude that ${\text{codim}}\, B(r', d')$ is very large when we assume $d$ to be conveniently large. \\

We seek now to establish an analogue of Lemma 1, that is to use the universal morphism
$\Phi: \o^N \rightarrow \v$ on $\bp \times C$ in order to give degeneracy loci representatives of the 
tautological classes \eqref{kun2}.\\

For each $i, \, 1 \leq i \leq r$, we choose a subbundle $\o^{r-i+1}$ of $\o^N$ on ${\bp} \times C.$ Denote by $\Delta_{i}$ and $\Delta_{i, p}$ the
noninjectivity loci of $\o^{r-i+1} \rightarrow \v$ on $\bp \times C$ and of $\o_p^{r-i+1} \rightarrow \v_p$ on $\bp \times \{p\}$ respectively.
$\Delta_{i}$ and $D_{i}$ agree of course on ${U} \times C$.  
We want to show as before

\begin{lemma}
(1) For a generic choice of trivial subbundle $\o^{r-i+1}$ of $\o^N$ on $\bp \times C$,
$$[\Delta_i] = c_i (\v) \cap [\bp \times C ],$$ $$[\Delta_{i,p}] = {\bar{a}}_i \cap [\bp]\, \,
\, \text{and}$$
$$[\eta (\Delta_i)] = {\bar{f}}_i \cap [\bp].$$
(2) Generic representative subschemes for all instances of ${\bar{a}}_i$, ${\bar{f}}_i$ and algebraic combinations of ${\bar{b}}_i^j$ classes in a monomial 
$P(\bar{a}, \bar{b}, 
\bar{f})$
intersect each other properly along $\mors \subset \bp$.
\label{deg222}
\end{lemma}

{\em{Proof.}} The proof is immediate. We need to show that $\Delta_{i}$ and $\Delta_{i,p}$ have
codimension $i$ in $\bp \times C$ and $\bp$ respectively, and that $\eta: \Delta_{i}
\rightarrow \eta(\Delta_{i})$ is generically one-to-one. But
${\left.\Delta_{i} \right |}_{U \times C}$ has codimension $i$ in ${U} \times C$
since it coincides with $D_{i}$ there. In the same way ${\left.\Delta_{i,p}
\right |}_{{U}}$ has codimension $i$ in ${U}$, and $\eta: \Delta_{i}
\rightarrow \eta(\Delta_{i})$ is generically one-to-one on ${U}$. On the other hand, 
the codimension of the boundary $\bp \setminus U = \cup_{r'< r\\, d'} B(r', d')$ is very large in $\bp$, hence the first part
of the lemma follows. The second part of the lemma follows from the second part of Lemma 1.\\

We aim to show that the intersection of generic degeneracy loci representatives of the
tautological classes appearing in a monomial $P (\bar{a}, \bar{b}, \bar{f}) {\bar{a}}_r^s,$ for $e = sr + t$, 
$t= \deg P$ and $t \leq \frac{N}{r}$
avoids the boundary $\sqcup_{r'< r, d'} B (r', d')$ of $\bp$. Proposition 1 ensures on the other hand that the intersection  
along $U \subset \quot$ is zero-dimensional. We will thus be able to conclude that the intersection
number $P(\bar{a}, \bar{b}, \bar{f}) {\bar{a}}_r^s \cap [\bp ]$ is given as the sum of multiplicities of the finitely many intersection points in $U$ 
of generic
degeneracy loci representatives of all the classes in $P(\bar{a}, \bar{b}, \bar{f}) {\bar{a}}_r^s$.\\

We analyze the intersection of a degeneracy locus
$\Delta_{r,p}$ with a boundary stratum $B (r', d')$. $\Delta_{r,p}$ is the zero locus of a section $\o \rightarrow \v_p$ on $\bp$.
Let $$0 \rightarrow {\e}'
\rightarrow \o^N \rightarrow {\f}' \rightarrow
0$$ be
the universal sequence on
$\text{Quot}_{d'} (\o^N, r', C) \times C.$
Let $\o$ be a subbundle of ${\o^N}^{\vee}$ (corresponding under $\Psi_{r', d'}$ to the one chosen before on $\bp$) 
and denote by $$\Delta_{p} (r', d') \subset \text{Quot}_{d'} (\o^N, r', C)$$ the
zero locus of the restricted homomorphism $\o_p \rightarrow {\e'}_{p}^{\vee}$
on
$\text{Quot}_{d'} (\o^N, r', C)) \times \{p\}.$ 
Given $\varphi: \o^N \rightarrow E,\;\;
\varphi \in B(r', d'),$ let $E' = {\text{Im}}\, \varphi.$ 
Note that $\o \rightarrow E'$ vanishes at $p$ if and only if $\o \rightarrow E$
does. Hence
$${\left. \Delta_{r,p} \right |}_{B(r', d')}\, = \,\Psi_{r', d'}^{-1} (\Delta_p
(r', d')).$$ We choose $g_k \in GL(N),\; p_k \in C,\;
1 \leq k \leq s$ so that
$g_k \Delta_{p_k} (r', d')$ intersect properly on $\text{Quot}_{d'} (\o^N, r',C).$
The intersection $$\cap_{k=1}^{s} g_k \Delta_{p_k} (r', d')$$
 is empty if its
codimension exceeds the dimension of the ambient scheme
$\text{Quot}_{d'} (\o^N, r', C).$ We can assume that $\text{Quot}_{d'} (\o^N,
r', C)$
has the expected dimension $N d' - r' (N-r') (g-1).$
The intersection on the other
hand has codimension 
$$s r' = \frac{e-t}{r} r' \geq \frac{Nd-r(N-r) (g-1) - \frac{N}{r}}{r} r'> Nd' - r(N-r') (g-1).$$
The last inequality is equivalent to 
\begin{equation}
\label{codim2}
N \left ( \frac{d}{r} - \frac{d'}{r'} \right ) > -(r-r') (g-1) + \frac{N}{r^2},
\end{equation}
which holds since $$N \left (\frac{d}{r} - \frac{d'}{r'} \right ) \geq \frac{N}{rr'}. $$ 

We conclude that the intersection of generic representative subschemes for the classes in the product $P(\bar{a},\bar{f},\bar{b})
{\bar{a}}_r^s$ avoids the boundary $\sqcup_{r'< r, d'} B (r', d')$ of $\bp$. The intersection is therefore contained in $U$, and since 
$U \subset \quot,$ by Proposition 1 we can assume that the intersection lies in fact in $\mors$, and is proper along $\mors.$
To summarize, we have obtained

\begin{proposition}  
The scheme-theoretic intersection of generic degeneracy loci representatives for the classes appearing in the top product $P(\bar{a},\bar{f},\bar{b}) 
{\bar{a}}_r^s, \, \deg P 
\leq \frac{N}{r},$ on
$\bp$ is zero-dimensional of length equal to the evaluation
$P(\bar{a},\bar{f},\bar{b}) {\bar{a}}_r^s \cap [\bp]$, and lies in $\mors \subset \bp.$
\end{proposition}

Since $\quot$ and $\bp$, as well as their universal structures, agree on $\mors$, Propositions 1 and 2 imply Theorem 1.

\end{document}